\documentclass[12pt]{amsart}
\usepackage{amssymb,amsmath,amsthm,amsfonts,amsopn,url,color,hyperref,enumerate,mathtools}
\usepackage[normalem]{ulem}
\usepackage{amscd}
\usepackage[mathscr]{euscript}
\usepackage{mathtools,rotating}
\usepackage{todonotes}
\usepackage{comment}
\usepackage[all,2cell]{xy}
\UseAllTwocells
\input xy
\xyoption{all}
\usepackage{pgf,tikz}
\usepackage{mathrsfs}
\usetikzlibrary{arrows}
\usetikzlibrary{patterns}
\usepackage{tikz,pgfplots}
\usetikzlibrary{matrix, cd}
\usepackage{esvect}
\usepackage{tabularx}
\usepackage{arydshln}
\usepackage{array}
\newtheorem{theorem}{Theorem}[section]
\newtheorem{lemma}[theorem]{Lemma}
\newtheorem{proposition}[theorem]{Proposition}

\newtheorem{corollary}[theorem]{Corollary}
\newtheorem{conjecture}[theorem]{Conjecture}
\theoremstyle{definition}\newtheorem{remark}[theorem]{Remark}
\newtheorem*{question*}{Question}
\theoremstyle{definition}\newtheorem{definition}[theorem]{Definition}
\newtheorem{example}[theorem]{Example}

\author{Hannah Klawa}
\address{Department of Mathematical Sciences \\ Indiana University East \\ Richmond, IN 47374}
\email{hklawa@iu.edu}

\author{Cigole Thomas}
\address{Department of Mathematics \\ Bates College \\ Lewiston, ME 04240}
\email{cthomas5@bates.edu}

\title[A Characterization of Completely Integrally Closed PVMDs]{A Characterization of Completely Integrally Closed Pr\"{u}fer $v$-Multiplication Domains}

\subjclass[2010]{Primary: 13F05, 13G05  Secondary: 13A15, 06F15}
\keywords{Archimedean domains, archimedean groups, completely integrally closed domains, lattice ordered groups, Pr\"{u}fer domains, PVMD, $t$-invertible ideals}

\begin{document}

\begin{abstract} We generalize a result known for Pr\"ufer domains to show that a PVMD $R$ is completely integrally closed if and only if $\text{Inv}^t(R)$ is an archimedean $\ell$-group.
It has been conjectured that a PVMD is completely integrally closed if and only if it is an
Archimedean domain. It is known that a completely integrally closed PVMD is an Archimedean domain.  We provide an example to illustrate that an Archimedean PVMD need not be completely integrally closed. 
\end{abstract}

\maketitle
\section{Introduction}

An integral domain $R$ is a \textit{Pr\"{u}fer domain} if every nonzero finitely generated ideal of $R$ is invertible. This characterization of Pr\"{u}fer domains and many other characterizations can be found in~\cite{Prufer} and~\cite{Gilmer}. Pr\"{u}fer $v$-multiplication domains, PVMDs, are a generalization of Pr\"{u}fer domains. An integral domain $R$ is a \textit{PVMD} if every nonzero finitely generated ideal of $R$ is $t$-invertible. For a more detailed discussion of PVMDs, the reader is referred to~\cite{Mott1981} and~\cite{Elliott19}.

An integral domain $R$ is \textit{Archimedean} if $\bigcap\limits_{n=1}^\infty (x^n)=(0)$ for all nonunits $x\in R$. Any completely integrally closed domain is Archimedean \cite[Corollary 13.4]{Gilmer}. But as noted in \cite{Anderson2017}, an Archimedean domain need not be completely integrally closed. In~\cite[Conjecture 3.2]{Anderson2017}, the authors made a conjecture (which they suspected could be false) that for PVMDs, being completely integrally closed is equivalent to being an Archimedean domain. 
For GCD domains (domains in which the $t$-closure of any non-zero finitely generated ideal is principal), being completely integrally closed is equivalent to being an Archimedean domain \cite[Theorem 3.1]{Dobbs1976}. In \cite{Anderson2017}, the authors give special classes of PVMDs where this conjecture holds which include $t$-QR PVMDs and integrally closed AGCD domains \cite[Corollary 3.11 and Corollary 3.12]{Anderson2017}. We provide an Example~\ref{counterexample} from~\cite[Example 5.9]{Gabelli2020} to show that \cite[Conjecture 3.2]{Anderson2017} is false: that is, an Archimedean PVMD need not be completely integrally closed. 

In~\cite[Proposition 3.1]{ICPrufer}, the authors show that a Pr\"{u}fer domain $R$ is completely integrally closed if and only if $\text{Inv}(R)$ is an archimedean $\ell$-group. We generalize this to show that a Pr\"{u}fer $v$-multiplication domain is completely integrally closed if and only if $\text{Inv}^t(R)$ is an archimedean $\ell$-group. 

\section{Preliminaries}
We begin by recalling some of the well-known definitions for the readers convenience which can be found in~\cite{Elliott19} or \cite{chapman2000}. Let $R$ be an integral domain with field of fractions $K$. Let $F(R)$ denote the set of all nonzero fractional ideals of $R$. The $v$-\textit{closure star operation} is defined by $I^v = (I^{-1})^{-1}$ for $I \in F(R)$. The \textit{$t$-closure star operation} is defined by 
$$I^t = \bigcup \left\{(J^v\, : \, J \in F(R) \text{ is finitely generated and } J \subseteq I\right\}.$$
For any $I \in F(R)$, $I^t \subseteq I^v$. An ideal $I\in F(R)$ is \textit{$t$-invertible} if $(II^{-1})^t=R$. 
An ideal $I$ is a \textit{$t$-ideal} if $I^t=I$.
Let $\text{Inv}^t(R)$ denote the set of all $t$-invertible $t$-ideals. Recall that $\text{Inv}^t(R)$ is a group under $t$-multiplication ($I\star_t J=(IJ)^t$) and admits a partial order given by $I\leq J$ if and only if $J\subseteq I$ for all $I, J \in \text{Inv}^t(R)$.

\begin{definition}
A group $G$ is \textit{archimedean} if whenever $a,b\in G$ with $na\leq b$ for all $n>0$, then $a\leq 0 $. 
\end{definition}
\begin{definition}
 An abelian group $(G,+)$ is a \textit{lattice ordered group} ($\ell-$ group) if it is partially ordered with a relation $\leq$ such that the meet $x\wedge y$ and join $x\vee y$ of $x$, $y\in G$ exist and such that for all $x, y, a \in G$, $x \leq y$ implies $a+x\leq a+y$. Throughout the paper, by an ``$\ell-$group" we mean an abelian $\ell-$group.
 \end{definition}
As in \cite{Zafrullah2000}, let $I \vee J = I\cap J$ and $I\wedge J=(I+J)^v$.

\begin{proposition}[\cite{Zafrullah2000}, Proposition 2.4]\label{Prop:l-group}
    An integral domain $R$ is a PVMD  if and only if $(\text{Inv}^t(R), \leq, \star_t) $ is a lattice-ordered abelian group. 
\end{proposition}
\begin{remark}\label{Rmk1}
Note that $\text{Inv}^t(R)$ is archimedean if whenever $I,J \in \text{Inv}^t(R)$ with $J\subseteq (I^n)^t$ for all $n>0$, then $R\subseteq I$.
\end{remark}
\begin{definition}
An integral domain $R$ is \textit{Archimedean} if $\bigcap\limits_{n=1}^\infty (x^n)=(0)$ for all nonunits $x\in R$. 
\end{definition}

\begin{definition}
An integral domain $R$ with quotient field $K$ is \textit{completely integrally
closed} if $ra^n \in R$ for $a \in K$, $0 \neq  r \in R$, and every integer $n \geq 1$ implies
$a \in R$. 
\end{definition}

\section{Completely integrally closed PVMDs and archimedean $\ell$-groups}

We will start by giving a generalization of \cite[Theorem 2.2 ($a\Rightarrow b)$]{Gilmer1972} for PVMDs below. This result might be well known but we include a proof for the reader's convenience.
\begin{lemma}\label{lemABn}
Let $R$ be an  PVMD. Then $((I \cap J)^n)^t = (I^n)^t \cap (J^n)^t$ where $I$ and $J$ are $t$-invertible $t$-ideals of $R$.

\end{lemma}

\begin{proof}
    Note that  $(I \cap J)^n = I^n \cap J^n$ holds for any valuation domain. Let $\{M_{\lambda}\}$ denote the set of $t$-maximal ideals of $R$. Because $R$ is a PVMD, $R=\mathop{\cap}\limits_{\lambda}R_{M_{\lambda}}$  where each $R_{M_{\lambda}}$ is a valuation ring \cite[Theorem 3.2.9, Corollary 3.2.12]{Elliott19}. It follows that
    \begin{eqnarray*}
    ((I \cap J)^n)^t & = & \mathop{\bigcap}\limits_{\lambda} (I \cap J)^n)R_{M_{\lambda}} \\
    & = & \mathop{\bigcap}\limits_{\lambda}[ (I \cap J)R_{M_{\lambda}}]^n \\
    &=& \mathop{\bigcap}\limits_{\lambda} (IR_{M_{\lambda}} \cap JR_{M_{\lambda}})^n \\
    &=& \mathop{\bigcap}\limits_{\lambda} [(IR_{M_{\lambda}})^n \cap (JR_{M_{\lambda}})^n] \\
    &=& \left[\mathop{\bigcap}\limits_{\lambda} (I^nR_{M_{\lambda}})\right] \bigcap \left[\mathop{\bigcap}\limits_{\lambda} (J^nR_{M_{\lambda}})\right] \\
    & = & (I^n)^t \cap (J^n)^t.
    \end{eqnarray*}
\end{proof}
We now give a standard result about $t$-invertible ideals in the lemma below. 

\begin{lemma}\label{Lemmat-v}
Let $R$ be an integral domain. Then
$\mathop{\cap}\limits_{n>0}(I^n)^t= \mathop{\cap}\limits_{n>0}(I^n)^v$ for proper $t$-invertible $t$-ideals. 
\end{lemma}

\begin{proof}
    By \cite[Proposition 1.2]{Zafrullah2000}, $I^n$ is $t$-invertible if $I$ is $t$-invertible. Because $t$-invertible $t$-ideals are divisiorial \cite[Theorem 1.1]{Zafrullah2000}, $I^t=I^v$ for all proper $t$-invertible $t$-ideals. 
    Consequently, $\mathop{\cap}\limits_{n>0}(I^n)^v = \mathop{\cap}\limits_{n>0}(I^n)^t$. 
\end{proof}

\begin{lemma}\label{Intersection}
    Let $R$ be a PVMD. Then $R$  is completely integrally closed if and only if $\bigcap\limits_{n>0}(I^n)^t=(0)$ for any proper $t$-invertible $t$-ideal.
\end{lemma}
\begin{proof}
The result follows by \cite[Theorem 1.2]{Anderson2017} and Lemma \ref{Lemmat-v}. 
\end{proof}
We prove the PVMD analog of the part of \cite[Proposition 3.1]{ICPrufer} which states that a Pr\"{u}fer domain $R$ is completely integrally closed if and only if $\text{Inv}(R)$ is an archimedean $\ell$-group.
\begin{proposition}\label{CIC-ALG}
A PVMD $R$ is completely integrally closed if and only
if $\text{Inv}^t(R)$ is an archimedean $\ell$-group. 
\end{proposition}

\begin{proof}
Suppose that $R$ is completely integrally closed with $I,J \in \text{Inv}^t(R)$ such that $J \subseteq (I^n)^t$
for all $n>0$. Then by Lemma \ref{lemABn}, we have $((I \cap R)^n)^t = (I^n)^t \cap R $. Hence $J \cap R \subseteq ((I \cap R)^n)^t = (I^n)^t \cap R$
for all $n>0$. 
Note that $J\cap R\neq (0)$ and so $\bigcap (I^n)^t \neq (0)$. 
Consequently, since $R$ is completely
integrally closed, it follows from Lemma \ref{Intersection} that $I\cap R=R$. Hence, $R\subseteq I$. Thus
$\text{Inv}^t(R)$ is an archimedean $\ell-$group by Remark \ref{Rmk1} and Proposition \ref{Prop:l-group}. 

Conversely, let $Inv^t(R)$ be an archimedean $\ell$-group and $I\in \text{Inv}^t(R)$ be a proper $t$-ideal. Then $\mathop{\cap}\limits_{n>0}(I^n)^t=(0)$ otherwise by Remark \ref{Rmk1}, $R\subseteq I$ since $\text{Inv}^t(R)$ is archimedean. 
Hence $R$ is completely
integrally closed by Lemma \ref{Intersection}. 
\end{proof}

One might wonder if for a PVMD $R$, there is a relation between being an Archimedean domain and $\text{Inv}^t(R)$ being an archimedean $\ell$-group. 
\begin{corollary}\label{Cor:l-group}
Let $R$ be a PVMD. If $\text{Inv}^t(R)$ is an archimedean $\ell$-group, then $R$ is an Archimedean domain. 
\end{corollary}
\begin{proof}
    Let $\text{Inv}^t(R)$ be an archimedean $\ell$-group. Then, by Proposition \ref{CIC-ALG}, $R$ is completely integrally closed. Hence, it follows by \cite[Corollary 3.4]{Anderson2017} that $R$ is Archimedean. 
\end{proof}

In~\cite{Anderson2017}, the authors define an integral domain $R$ to be \emph{strongly Archimedean} if $\mathop{\cap}\limits_{n =1}^\infty (a/b)^n = (0)$ for every $a, b \in R$ with $(b) \not\subseteq (a)$. In~\cite[Proposition 3.3]{Anderson2017}, they show that for any integral domain being strongly Archimedean is equivalent to being completely integrally closed. So, we obtain the following corollary of Proposition~\ref{CIC-ALG}.

\begin{corollary}
Let $R$ be a PVMD. Then $R$ is strongly Archimedean if and only if $\text{Inv}^t(R)$ is an Archimedean $\ell$-group. 
\end{corollary}
\begin{proof}
This follows from Proposition \ref{CIC-ALG} and \cite[Proposition 3.3]{Anderson2017}.
\end{proof}
\section{Archimedean PVMDs need not be completely integrally closed}
In \cite{Anderson2017}, the authors made the following conjecture which they suspected could be false. 
\begin{conjecture}[\cite{Anderson2017}, Conjecture 3.2]
    Let $R$ be a PVMD. Then $R$ is completely integrally closed if and only if $R$ is Archimedean.
\end{conjecture}

The following example shows that this conjecture is false and that the converse of Corollary \ref{Cor:l-group} does not hold. 
\begin{example}\label{counterexample}
Consider the example by Gabelli and Roitman from \cite[Example 5.9]{Gabelli2020}. 
Let $X$ and $Y$ be two independent indeterminates over a field $k$. Let
$$R=k\left[X,Y,  \frac{X(1-X)^n}{Y^n}, \  \frac{Y^{n+1}}{(1-X)^n}\ (n \in \mathbb{N})\right]_S,$$
where $S=k[Y] \backslash Yk[Y]$.
Let $T=\frac{1-X}Y$, 
so
$$R=k\left[Y, YT, (1-YT)T^n, \frac Y{T^n}\ (n\in \mathbb{N})\right]_S.$$
By \cite[Example 5.9]{Gabelli2020}, $R$ is an Archimedean Pr\"{u}fer domain and hence an Archimedean PVMD. 
As noted in \cite[Example 5.9]{Gabelli2020}, $T\notin R$. To see why, suppose $T$ is in $R$. Note that $R \subseteq k[Y, T, T^{-1}]_S$. Let 
$$\varphi: k[Y, T, T^{-1}]_S \longrightarrow k[T, T^{-1}]$$
be given by $$Y \mapsto 0, \ T \mapsto T \ \text{ and } \ r \mapsto r$$ for all $r \in k$. Let $h/f = h^\prime/g$ where $f, g \in S$ and $h, h^\prime \in k[Y, T, T^{-1}]$. Then $f, g \in k[Y]$ such that $f(0) \neq 0, g(0) \neq 0 \in k$. It follows that $\varphi(h/f) =\varphi(h'/g) \in  k[T, T^{-1}]$ and hence $\varphi$ is well-defined. Note that  $$Y \mapsto 0, \ (1-YT)T^n \mapsto T^n, \ Y/T^n \mapsto 0$$ and $\varphi(S) \subseteq k\backslash \{0\}$. Hence, $\varphi(R) \subseteq k[T]$. Because $1/T \in k[Y, T, T^{-1}]_S$ and $\varphi(1/T) = 1/T \not\in k[T]$, it follows that $1/T \not \in R$. So, $T$ is not a unit in $R$. Observe that $Y = T^n\frac{Y}{T^n}$ and $\frac{Y}{T^n} \in R$ and so $Y \in (T^n)$ for all $n$ which implies $Y \in \bigcap\limits_{n=1}^\infty (T^n)$. 
Because $T$ is not a unit, it follows that $Y = 0$ since $R$ is Archimedean which contradicts that $Y \neq 0 \in R$. Thus, $T \not\in R$.

Note that $XT^n\in R$ for all $n\in \mathbb{N}$ but $T \not\in R$. It follows that $R$ is not completely integrally closed. 
Thus, an Archimedean PVMD need not be completely integrally closed. 
\end{example}

\begin{corollary}
    Let $R$ be a PVMD. Then $R$ being Archimedean does not imply that $\text{Inv}^t(R)$ is an archimedean $\ell$-group.
\end{corollary}
\begin{proof}
By Proposition \ref{CIC-ALG}, $\text{Inv}^t(R)$ is not an archimedean $\ell$-group for the domain $R$ from Example $\ref{counterexample}$. This implies that the converse of Corollary \ref{Cor:l-group} is false. 
\end{proof}
\section*{Acknowledgments}

The first named author is grateful for the Faculty Research and Creative Activity Support Fund at Indiana University East which helped provide support during part of this study. We would also like to thank Neil Epstein for a useful discussion about Example~\ref{counterexample}.

\bibliographystyle{plain}
\bibliography{references}

@Book{Gilmer,
author = {Gilmer, R.},
title = {Multiplicative Ideal Theory},
publisher = {Marcel Decker, Inc.},
year = {1972}
}

@Book{Prufer,
author = {Marco Fontana and James A. Huckaba and Ira J. Papick},
title = {Pr\"{u}fer Domains},
publisher = {Marcel Dekker, Inc.},
year = {1997}
}

@article {ICPrufer,
    AUTHOR = {Heubo-Kwegna, Olivier A. and Olberding, Bruce and Reinhart,
              Andreas},
     TITLE = {Group-theoretic and topological invariants of completely
              integrally closed {P}r\"ufer domains},
   JOURNAL = {J. Pure Appl. Algebra},
  FJOURNAL = {Journal of Pure and Applied Algebra},
    VOLUME = {220},
      YEAR = {2016},
    NUMBER = {12},
     PAGES = {3927--3947},
      ISSN = {0022-4049,1873-1376},
   MRCLASS = {13F05 (06F15 13A05 13A15)},
  MRNUMBER = {3517563},
MRREVIEWER = {Daniel\ D.\ Anderson},
       DOI = {10.1016/j.jpaa.2016.05.021},
       URL = {https://doi.org/10.1016/j.jpaa.2016.05.021},
}

@book{Elliott19,
  author    = {Jesse Elliott},
  title     = {Rings, Modules, and Closure Operations},
  series    = {Springer Monographs in Mathematics},
  publisher = {Springer, Cham},
  year      = {2019},
  isbn      = {978-3-030-24400-2},
  eisbn     = {978-3-030-24401-9},
  doi       = {10.1007/978-3-030-24401-9}
}

@Inbook{Zafrullah2000,
author="Zafrullah, Muhammad",
editor="Chapman, Scott T.
and Glaz, Sarah",
title="Putting T-Invertibility to Use",
bookTitle="Non-Noetherian Commutative Ring Theory",
year="2000",
publisher="Springer US",
address="Boston, MA",
pages="429--457",
isbn="978-1-4757-3180-4",
doi="10.1007/978-1-4757-3180-4_20",
url="https://doi.org/10.1007/978-1-4757-3180-4_20"
}

@article{Anderson2017,
author = {D. D. Anderson and David F. Anderson and Muhammad Zafrullah},
title = {Completely integrally closed {P}rüfer v-multiplication domains},
journal = {Communications in Algebra},
volume = {45},
number = {12},
pages = {5264--5282},
year = {2017},
publisher = {Taylor \& Francis},
doi = {10.1080/00927872.2017.1303502},


URL = { 
    
        https://doi.org/10.1080/00927872.2017.1303502
    
    

},
eprint = { 
    
        https://doi.org/10.1080/00927872.2017.1303502
    
    

}

}

@article{Gabelli2020,
author = {Stefania Gabelli and Moshe Roitman},
title = {{On finitely stable domains, II}},
volume = {12},
journal = {Journal of Commutative Algebra},
number = {2},
publisher = {Rocky Mountain Mathematics Consortium},
pages = {179 -- 198},
year = {2020},
doi = {10.1216/jca.2020.12.179},
URL = {https://doi.org/10.1216/jca.2020.12.179}
}

@article {Gilmer1972,
    AUTHOR = {Gilmer, Robert and Grams, Anne},
     TITLE = {The equality {$(A\cap B)\sp{n}=A\sp{n}\cap B\sp{n}$} for
              ideals},
   JOURNAL = {Canadian J. Math.},
  FJOURNAL = {Canadian Journal of Mathematics. Journal Canadien de
              Math\'ematiques},
    VOLUME = {24},
      YEAR = {1972},
     PAGES = {792--798},
      ISSN = {0008-414X,1496-4279},
   MRCLASS = {13G05},
  MRNUMBER = {308109},
MRREVIEWER = {J.\ W.\ Brewer},
       DOI = {10.4153/CJM-1972-076-9},
       URL = {https://doi.org/10.4153/CJM-1972-076-9},
}

@book{chapman2000,
  title={Non-{N}oetherian Commutative Ring Theory},
  editor={Chapman, Scott T. and Glaz, Sarah},
  series={Mathematics and Its Applications},
  volume={520},
  year={2000},
  publisher={Springer New York, NY},
  doi={https://doi.org/10.1007/978-1-4757-3180-4}
}

@article{Mott1981,
  author  = {Mott, Joe L. and Zafrullah, M.},
  title   = {On Pr{\"u}fer $v$-multiplication domains},
  journal = {Manuscripta Mathematica},
  volume  = {35},
  year    = {1981},
  pages   = {1-26},
  doi     ={https://doi.org/10.1007/BF01168446}
}

@article{Dobbs1976,
title={On a Class of Archimedean Integral Domains}, volume={28}, DOI={10.4153/CJM-1976-038-x}, number={2}, journal={Canadian Journal of Mathematics}, author={Beauregard, Raymond A. and Dobbs, David E.}, year={1976}, pages={365–375}}

\end{document}